\documentclass{amsart}
\usepackage{euscript,eufrak,amscd,epsfig,verbatim}

\newcommand{\sh}[1]{\EuScript{#1}}     
\newcommand{\goth}[1]{\EuFrak{#1}}   
\newcommand{\mb}[1]{\mathbf{#1}}       
\newcommand{\C}{\mb{C}}                
\newcommand{\Q}{\mb{Q}}                
\newcommand{\Z}{\mb{Z}}                
\newcommand{\ctop}{c_{\mathrm{top}}}   
\newcommand{\isoa}{\xra{\sim}}         
\renewcommand{\P}[1][{}]{\ensuremath{\mb{P}^{#1}}}
\DeclareMathOperator{\Sym}{Sym}
\DeclareMathOperator{\Hom}{Hom}
\renewcommand{\H}{\operatorname{H}}
\DeclareMathOperator{\codim}{codim}
\DeclareMathOperator{\CH}{CH}
\DeclareMathOperator{\ev}{ev}
\DeclareMathOperator{\id}{id}
\newcommand{\hra}{\hookrightarrow}
\newcommand{\iso}{\cong}

\newcommand{\ol}[1]{\overline{#1}}

\newcommand{\xra}[1]{\xrightarrow{#1}}
\renewcommand{\epsilon}{\varepsilon}
\newcommand{\Gbar}{\ol{G}}
\newcommand{\Mbar}{\ol{M}}

\newcommand{\M}[1]{\Mbar_{0,#1}}
\newcommand{\G}[1]{\Gbar_{0,#1}}
\newcommand{\T}{\mb{T}}
\newcommand{\vir}[1]{\left[#1\right]^{\mathrm{vir}}}
\theoremstyle{plain}
\newtheorem{thm}{Theorem}[section]

\newtheorem{lemma}[thm]{Lemma}
\newtheorem{cor}[thm]{Corollary}
\newtheorem{formula}[thm]{Formula}
\theoremstyle{definition}
\newtheorem{defn}{Definition}
\newtheorem{exas}{Examples} 
\newtheorem{nota}{Notation} 
\newtheorem{ISC}{Important Special Case}
  
\newtheorem{observation}{Observation}

\theoremstyle{remark}
\newtheorem{remark}{Remark}  
\newtheorem{claim}{Claim} 


\title{New recursions for genus-zero {G}romov-{W}itten invariants}
\author{Aaron Bertram}
\address{Dept. of Mathematics, University of Utah, 155 S. 1400
E., Salt Lake City, UT 84112}
\email{bertram@math.utah.edu}
\thanks{The first author was supported in part by NSF Research Grant
DMS-9970412} 
\author{Holger P. Kley}
\address{Dept. of Mathematics, Colorado State University, Fort
Collins, CO 80523}
\email{kley@math.colostate.edu}


\begin{document}

\addtocounter{section}{-1}

\begin{abstract}
New relations among the genus-zero Gromov-Witten invariants of a
complex projective manifold $X$ are exhibited. When the cohomology of
$X$ is generated by divisor classes and classes ``with vanishing
one-point invariants,'' the relations determine many-point invariants in
terms of one-point invariants. 
\end{abstract}

\maketitle

\section{Introduction} 
The localization theorem for equivariant cohomology has recently been
used with great success  to compute the genus-zero Gromov-Witten
invariants relevant to the mirror conjecture \cite{Giv,LLY,B1}.
Zero-point invariants count expected numbers 
of rational curves on a projective manifold~$X$, while 
the more general $m$-point invariants count expected numbers of
rational curves meeting $m$ given submanifolds (or 
cohomology classes).  For the the mirror conjecture, only
the zero and one-point invariants are computed, though for the
construction of the quantum product (even the
small version), one needs more general invariants.

\medskip

In this paper we will apply the localization theorem to  
study genus-zero Gromov-Witten invariants involving any number
of marked points. A straightforward generalization of Givental's
(one variable) $J$-function yields homology-valued $J$-functions
in any number of variables $t_1$,\dots,$t_{m}$ which encode all the
(generalized) $m$-point genus-zero invariants.
Our main theorem is a collection of relations among these
$J$-functions expressing a part of the $J$-function for a fixed
curve class and number of variables in terms of the $J$-functions
involving fewer variables and/or ``smaller'' curve classes.
When the cohomology of $X$ is generated by divisor classes,
or, more generally, when every class orthogonal to the subring
generated by divisor classes annihilates (via cap product) all
one-variable $J$-functions, then these new
relations  completely determine all $m$-point genus-zero Gromov-Witten
invariants (of classes generated by divisor classes) in terms of
one-point invariants. That is, in this setting, 
the one-variable $J$-function determines all the 
others.  A complete intersection in
$\P[n]$ has this orthogonality property, and in that case we
exhibit a formula expressing ``mixed'' two-point
invariants in terms of one-point invariants. We apply this new
formula to compute previously unknown quantum products of cohomology
classes on  Fano complete intersections, where the one-variable
$J$-function is known. Since our recursions do not require any
positivity of $X$, they would apply just as well to general-type
complete intersections. Unfortunately, in those cases, the
one-variable $J$-function is not known.

\medskip

The idea is the following. Given a stable map with several marked
points, a copy of $\P[1]$ is attached to each of the marked
points (at $0$ on  the $\P[1]$).
This allows one to embed the moduli space of 
$m$-pointed stable maps into the moduli space
of stable maps with no marked points and $m$ 
parametrizations.  We will call the latter space the {\sl graph
space.}  There is a natural torus action on the
graph space, one of whose fixed loci is the given
moduli space of stable maps. There are
equivariant forgetful morphisms among the graph spaces, and by 
comparing residues of some carefully chosen equivariant cohomology
classes along the fixed loci, we  obtain the recursive
formulas for the $J$-functions.  A startling (to us) feature of this
approach is that it is much simpler than the computation of
one-point invariants, since our argument requires no analysis of the
boundary of graph spaces. 

\section{Kontsevich-Manin spaces} 
We recall the basic
properties of the genus-zero stable map spaces and 
some results on Gromov-Witten invariants, and give an
instance of our formula (to be proved in greater generality later).

\begin{defn} 
A morphism $f\colon(C;p_1,\dots,p_m) \to X$ from a connected
$m$-pointed complex rational curve $C$ to a 
complex projective manifold $X$ is {\sl prestable} if $C$ has only
nodes as singularities and  $p_1$,\dots,$p_m \in C$ are nonsingular.
If in addition every irreducible component of $C$ collapsed by $f$ has 
three or more distinguished points---a {\sl distinguished point} is a node
or marked point---we say that $f$ is {\sl stable.} 
\end{defn}

\begin{remark} This notion of stability is analogous to
Deligne-Mumford stability for pointed curves. Indeed, a stable map
to a point is a stable pointed curve. 
\end{remark}

The moduli of stable maps has been extensively studied,
ever since stable maps were introduced by Kontsevich and
Manin~\cite{KM}.  (See also \cite{BM} and \cite{FP} as good references
for the following properties.)

\medskip

Given $\beta \in \H_2(X,\Z)$, there is a proper Deligne-Mumford
stack $\M{m}(X,\beta)$, representing
flat families of genus-zero stable maps with $m$ marked points and
image homology class~$\beta$. 
Each of the moduli stacks
$\M{m}(\P[n],d)$ is smooth (as a stack)
of the expected dimension
$n + (n+1)d + (m-3)$. For general $X$, there is a virtual
class
$\vir{\M{m}(X,\beta)}$ in the Chow group of 
$\M{m}(X,\beta)$ of the expected dimension
$n - \deg_{K_X}(\beta) + (m-3)$. 

\medskip

There are forgetful maps and evaluation maps:
\begin{equation*}\begin{CD}
  \M{m+1}(X,\beta) @>{e_i}>> X \\
  @V{\pi_i}VV @.\\
  \M{m}(X,\beta) @. {}
\end{CD}\end{equation*}
where $\pi_i$ ``forgets'' the marked point $p_i$ (and collapses
components, if necessary), and
$e_i$ evaluates the stable map at~$p_i$. When
$i = m+1$, this diagram can be taken as part of the ``universal
stable map'' over~$\M{m}(X,\beta)$. The rest of the
universal stable map consists of sections
\[ 
  \rho_i\colon \M{m}(X,\beta)\to \M{m+1}(X,\beta)
\]
of $\pi_{m+1}$ corresponding to the marked points. 

\medskip

In case
$X \subset \P[n]$ is the transverse zero locus of a section of a vector
bundle $E$ on $\P[n]$ which is generated by global sections,
the refined top Chern class 
$\ctop({\pi_{m+1}}_*e_{m+1}^*E)$ on $\M{m}(\P[n],d)$ produces the
virtual class on~$\M{m}(X,d)$.

\medskip

Morphisms
$\phi\colon X \to Y$ give rise to morphisms of stable map
spaces
\begin{equation*}
  \M{m}(X,\beta) \to \M{m}(Y,\phi_*\beta).
\end{equation*}
Finally, the ``boundary'' of $\M{m}(X,\beta)$ is covered by
the images of the gluing maps:
\begin{equation*}
  \delta_{S,\alpha} \colon \M{k+1}(X,\alpha) \times_X
  \M{m-k+1}(X,\beta - \alpha) \to \M{m}(X,\beta)
\end{equation*}
where $S\subseteq \{1,\dots,m\}$ is a subset of cardinality~$k$,
which,  together with the curve class~$\alpha$,  
describes how the stable map breaks into (at least) two components. 

\medskip

The {\sl Gromov-Witten invariants} are usually
interpreted as intersection numbers on the  Kontsevich-Manin spaces
of stable maps. Given cohomology classes
$\gamma_1$,\dots,$\gamma_m$ on~$X$, one defines the ``ordinary''
invariants
\begin{equation*}
  \left<\gamma_1,\dots,\gamma_m\right>^X_\beta  :=
  \deg\left(\pi_1^*\gamma_1 \cup \dotsb \cup
  \pi_m^*\gamma_m \cap \ev_*\vir{\M{m}(X,\beta)}\right),
\end{equation*}
where 
\[
\ev := (e_1,\dots,e_m) \colon \M{m}(X,\beta) \to X^m
\] 
is the total evaluation map and $\pi_i \colon X^m\to X$ 
are the projections. 

The general invariants are defined using the {\sl cotangent classes}
\begin{equation*}
  \psi_i := c_1(\rho_i^*\omega_{\pi_{m+1}}),
\end{equation*}
where $\omega_{\pi_{m+1}}$ is the relative dualizing sheaf. The
general invariants are:
\begin{multline*}
  \left<\gamma_1\psi^{a_1},\dots,\gamma_m\psi^{a_m}\right>^X_\beta := \\
  \deg\left(\pi_1^*\gamma_1 \cup \dotsb \cup
  \pi_m^*\gamma_m \cap \ev_*(\psi_1^{a_1}
  \cup \dotsb \cup \psi_m^{a_m} \cap \vir{\M{m}(X,\beta)})\right),
\end{multline*}
where $a_1$,\dots,$a_m$ are non-negative integers.

\medskip

The following is a very useful way to package $2$-point
``mixed'' invariants:
\begin{equation*}
  \left< \gamma_1,\frac{\gamma_2}{t-\psi}\right>^X_\beta :=
  t^{-1}\left< \gamma_1,\gamma_2\right>^X_\beta + t^{-2}\left<
  \gamma_1,\gamma_2\psi\right>^X_\beta +
  t^{-3}\left<\gamma_1,\gamma_2\psi^2\right>^X_\beta +\dotsb,
\end{equation*}
where $t$ is a variable.  Similarly, for general $1$-point invariants:
\begin{equation*}
  \left< \frac\gamma{t(t- \psi)}\right>^X_\beta := 
  t^{-2}\left< \gamma\right>^X_\beta + t^{-3}\left<
  \gamma\psi\right>^X_\beta + t^{-4}\left<
  \gamma\psi^2\right>^X_\beta +\dotsb .
\end{equation*}
The simplest of our formulas is expressed in terms of these packages
(extended $t$-linearly):

\begin{formula}\label{form:one}
Suppose $X\subset \P[n]$ is a complete
intersection of dimension $r \ge 3$ and degree~$l$. Then for  all
$0 \le a,b\le m$ and $d \ge 0$,
\begin{multline*}
  \left<H^a,\frac{H^b}{t-\psi}\right>^X_d +
    \left<\frac{H^a(H-dt)^b}{-t(-t-\psi)}\right>^X_d +\\
  \sum_{e=1}^{d-1}\sum_{c=0}^r\frac1{l} 
    \left<H^a,\frac{H^c}{t-\psi}\right>^X_{d-e}
    \left<\frac{H^{r-c}(H-et)^b}{-t(-t-\psi)}\right>^X_e \in \Q[t]
\end{multline*}.
\end{formula}
This formula implies a special case of our 
reconstruction theorem~\ref{thm:recon}:

\begin{cor} The mixed
two-point invariants of complete intersections in $\P[n]$ involving
only powers of $H$ are determined by the one-point invariants.
\end{cor}

\begin{proof} The first term in the formula is clearly determined
by the others. By induction on~$d$, the mixed
invariants of degree $d$ are therefore determined by the one-point
invariants of degree
$d$ or less.
\end{proof}

In the appendix, we use Formula~\ref{form:one} to compute small quantum
products of cohomology classes on Fano complete intersections. For
now, we point out the identities that follow from the formula
when the classes (in the  second slot) are of
codimensions $0$, $1$ and~$2$.

\begin{description}
\item[$\codim 0$] $\left<H^a,1\right>^X_d = 0$. 

\item[$\codim 1$] $\left<H^a,\psi\right>^X_1 = -\left<H^a\right>^X_d$
and $\left<H^a,H\right>^X_d = d\left<H^a\right>^X_d$.

\item[$\codim 2$]
\begin{align*}
  \left<H^a,\psi^2\right>^X_d & = \left<H^a\psi\right>^X_d -
    \frac1{l}\sum_{e=1}^{d-1}\sum_{c=0}^r 
    \left<H^a,H^c\right>^X_{d-e}\left<H^{r-c}\right>^X_e \\
  \left<H^a,H\psi\right>^X_d & = -\left<H^{a+1}\right>^X_d - 
    d\left<H^a\psi\right>^X_d + \frac1{l} \sum_{e=1}^{d-1}\sum_{c=0}^r
    e\left<H^a,H^c\right>^X_{d-e}\left<H^{r-c}\right>^X_e \\
  \left<H^a,H^2\right>^X_d & = 2d\left<H^{a+1}\right>^X_d +
    d^2\left<H^a\psi\right>^X_d - \frac1{l} \sum_{e=1}^{d-1}\sum_{c=0}^r
    e^2\left<H^a,H^c\right>^X_{d-e}\left<H^{r-c}\right>^X_e.
\end{align*}
\end{description}

Notice that the codimension two identities are not self-contained,
since they inductively involve classes of higher codimension. 

\begin{remark}
The codimension $0$ and the two codimension $1$ identities are special
cases of the string, dilaton and divisor equations, respectively.
The identities for codimension~$2$ classes, however,
are not special cases of any general equations that we are aware of
(though we've been informed that Lee and Pandharipande
\cite{LP} have another method for producing such identities).  
\end{remark}

To state our main theorem, we will use
$J$-functions of several variables, generalizing Givental's
one-variable definition.  
\begin{defn}
\begin{align*}
  J^X_\beta(t_1,\dots,t_m) &:= 
   \ev_*\left(\frac{\vir{\M{m}(X,\beta)}}{t_1(t_1-\psi_1)\dotsm
   t_m(t_m-\psi_m)}\right) \\
  & := \ev_*\left(\prod_{i=1}^mt_i^{-2} \left(1 + \frac{\psi_i}{t_i} +
   \frac{\psi_i^2}{t_i^2} + \dotsb\right)\cap\vir{\M{m}(X,\beta)}
   \right) \\
  & \in \H_*(X^m,\Q)[t_1^{-1},\dots,t_m^{-1}]
\end{align*}
with initial conditions:
\begin{equation*}
  J^X_0(t_1) := [X]
\end{equation*}
and
\begin{equation*}
  J^X_0(t_1,t_2) := \frac\Delta{t_1t_2(t_1+t_2)},
\end{equation*}
where $\Delta \in \H_*(X\times X,\Q)$ is the  diagonal
class.   
\end{defn}
\begin{remark}
  When $m=1$, our $J$-function is the Poincar\'e dual of Givental's.
\end{remark}

The $J$-functions encode all genus-zero Gromov-Witten
invariants.  The following result concerning the one-variable
$J$-function was first proved in \cite{Giv};  see \cite{LLY}
and \cite{B1} for alternate approaches.
\begin{thm}[Givental]\label{thm:giv}
\begin{enumerate}
\item If $X \subset \P[n]$ is a complete intersection of
type $(l_1,\dots,l_m)$ which is Fano of index two or more (i.e.
$l_1+\dots+l_m < n$), let $H$ be the hyperplane class. Then
\begin{equation*}
  J^X_d(t) = I^X_d(t) :=
  \frac{\prod_{i=1}^m\prod_{k=1}^{dl_i}(l_iH +
  kt)}{\prod_{k=1}^d(H+kt)^{n+1}}\cap [X].
\end{equation*}
\item If $l_1+\dots + l_m = n$ or $n+1$,  
then the following generating functions
coincide after an explicit  ``mirror
transformation'' (see \cite{Giv}, \cite{LLY} or
\cite{B1}):
\begin{equation*}
  J^X(q) := \sum_{d=0}^\infty J^X_d(t)q^d\quad\text{and}\quad I^X(q) := 
  \sum_{d=0}^\infty I^X_d(t)q^d.
\end{equation*}
\end{enumerate}
\end{thm}

We introduce a tool for manipulating $J$-functions of several variables: 
\begin{defn} Given classes $\Gamma_1
\in \H_*(X^k\times X,\Q)$, $\Gamma_2 \in \H_*(X \times X^{m-k},\Q)$
and $\gamma\in \H^*(X,\Q)$, we use the K\"unneth formula and 
Poincar\'e duality to regard the
tensor product as a $\Q$-linear map:
\begin{equation*}
  \Gamma_1\otimes\Gamma_2\colon \H^*(X^2,\Q)\to \H_*(X^m,\Q)
\end{equation*}
where the factors of $X^2$ are the distinguished 
factors of $X^k \times X$ and~$X\times X^{m-k}$. We then define
the {\sl twisted product}
$\Gamma_1\otimes_\gamma\Gamma_2
\in \H_*(X^m,\Q)$ by setting:
\begin{equation*}
  \Gamma_1\otimes_\gamma \Gamma_2 := (\Gamma_1 \otimes
  \Gamma_2)(\delta \cup \pi^*\gamma)
\end{equation*} 
where $\delta$ is the diagonal class and $\pi\colon X^2\rightarrow X$
is either of the two projections.
\end{defn}

\begin{exas}  Let $\gamma_1$, $\gamma_2$ be Poincar\'e dual to
$\Gamma_1$, $\Gamma_2$. 
\begin{enumerate} 
\item If $k = m = 0$, then $\Gamma_1 \otimes_\gamma 
\Gamma_2 \in \Q$ is the triple intersection:
\begin{equation*}
\int_X \gamma \cup \gamma_1 \cup \gamma_2.
\end{equation*}
More generally, if $k = m$, then $\Gamma_1 \otimes_\gamma \Gamma_2
= \pi^X_*(\pi_X^*(\gamma \cup \gamma_2)\cap \Gamma_1)$, where
$\pi_X\colon X^m \times X \to X$ and
$\pi^X\colon X^m \times X\to X^m$ are the two projections.

\item If $k = m-1$, then $\Gamma_1\otimes_\gamma \Delta = 
\pi_X^*\gamma\cap\Gamma_1$ for $\pi_X\colon X^{m-1}\times X
\to X$.
\end{enumerate}
\end{exas}

\begin{thm}[The main theorem---rank one case]\label{thm:main1} If $H$ is
an ample divisor class generating $\H^2(X,\Q)$ on a complex
projective manifold~$X$, then for each choice of
$m > 0$, $d = \deg_H(\beta) \ge 0$ and
$0 \le b \le \dim(X)$, 
\begin{multline*}
  t(t_1+t)\bigg( \sum_{1\in S\subseteq [m]} \sum_{e=0}^d
   J^X_{d-e}(\vec t_S,t)\otimes_{(H-et)^b}
   J^X_e(-t,\vec t_{S^c}) + \\
  \sum_{j=2}^m
   J^X_d(\vec t_{\hat j},t_j)\otimes_{H^b}
   J^X_0(-t_j,t)\bigg) \in \H_*(X^m,\Q)[t_1,t_1^{-1},\dots,t_m,t_m^{-1},t].
\end{multline*}
\end{thm}

\begin{nota} We set $[m] := \{1,\dots,m\}$. For 
subsets $S = \{s_1,\dots,s_k\} \subseteq [m]$, we define $\vec
t_S := (t_{s_1},\dots,t_{s_k})$. Since the $J$-functions are symmetric
in their variables, the expressions $J^X_d(\vec t_S,t)$ are
well-defined. We also set $\hat j := [m] - \{j\}$. 
\end{nota}

We will prove the main theorem later, as well as a more general
version where the rank one condition on $\H^2(X,\Q)$ is
removed. To finish this section, we show how Formula~\ref{form:one} follows
from the main theorem.

\begin{proof}[Proof of Formula~\ref{form:one}] We apply the main
theorem in the case $m = 1$. Only the double sum appears, in this case
as the single sum
\begin{equation}\label{eq:star}
  \sum_{e=0}^d J^X_{d-e}(t_1,t)\otimes_{(H - et)^b}J^X_e(-t).
\end{equation} 
The first and last of the terms are:
\begin{equation*}
  J^X_d(t_1,t)\otimes_{H^b}J^X_0(-t) = 
  J^X_d(t_1,t)\otimes_{H^b}[X] = {\pi_1}_*(J^X_d(t_1,t)
  \cup \pi_2^*H^b)
\end{equation*}
and
\begin{equation*}
  J^X_0(t_1,t)\otimes_{(H-dt)^b}J^X_d(-t) =
  \frac{\Delta \otimes_{(H-dt)^b} J^X_d(-t)}{t_1t(t_1+t)} = 
  \frac {(H-dt)^b \cup J^X_d(-t) }{t_1t(t+t_1)}.
\end{equation*}
Multiply \eqref{eq:star} through by $t_1t(t_1+t)$, and the main theorem
tells us we obtain an element of $\Q[t_1,t_1^{-1},t]$
when we integrate against $H^a$. For example, by the
projection formula, the first term gives
\begin{equation*}
  t_1t(t_1+t)\deg\big( H^a \cap {\pi_1}_*(\pi_2^*H^b\cap J^X_d(t_1,t))\big) = 
  (t_1+t)\left<\frac{H^a}{t_1-\psi},\frac{H^b}{t-\psi}\right>^X_d
\end{equation*}
and similarly for the other the terms.  When we consider only
the terms that are constant in~$t_1$, we obtain the following
formula:
\begin{multline*}
  \left<H^a,\frac{H^b}{t-\psi}\right>^X_d +
    \left<\frac{H^a(H-dt)^b}{-t(-t-\psi)}\right>^X_d +\\
  \sum_{e=1}^{d-1}\sum_{i,j=1}^N 
    \left<H^a,\frac{\gamma_i}{t-\psi}\right>^X_{d-e}g^{ij}
    \left<\frac{\gamma_j(H-et)^b}{-t(-t-\psi)}\right>^X_e \in \Q[t],
\end{multline*}
where $\gamma_1,\dots,\gamma_N \in \H^*(X,\Q)$ are a basis,
with respect to which $g^{ij}$ is the inverse of the intersection
matrix. This much holds for any ample $H$ generating $\H^2(X,\Q)$.

\medskip

The fact that $X$ is a complete intersection tells us that all the
one-point invariants of the form
$\left<\gamma\psi^c\right>^X_\beta$ vanish when $\gamma$ is a
primitive cohomology class. This can either be seen using 
Givental's formulas, or by a monodromy argument. Since a
basis for the cohomology may be chosen consisting of powers of $H$
and (orthogonal) primitive classes, this tells us that we may
replace the basis $\{\gamma_i\}$ by the smaller set $\{H^c\}$ of powers
of $H$, resulting in Formula~\ref{form:one}.
\end{proof}

\section{Graph spaces}  
Graph spaces are particular Kontsevich-Manin
spaces which come equipped with a natural torus
action. In this section, will describe some of the fixed
loci under this torus action in order to eventually apply 
the Atiyah-Bott localization theorem to prove the main theorem.

\begin{defn} The $m${\sl-parametrized graph space} is
\begin{equation*}
  \G{m}(X,\beta) := \M{0}(X\times (\P[1])^m,(\beta,1^m)).
\end{equation*}
\end{defn}

It is often useful to think of the graph space in the following
way. Let a {\sl parametrization} of a rational curve
$C$ be an isomorphism from $\P[1]$ to one of the 
irreducible components of~$C$.  A morphism $f\colon
(C;\P[1]_1,\dots,\P[1]_m) \to X$ from a connected rational
curve with $m$ parametrizations is{\sl  prestable} if $C$ has only
nodes as singular points. If in addition every component of $C$ is
either parametrized (possibly in several ways) or has at least three
nodes (or both), we say the $f$ is {\sl stable.} 

Then $\G{m}(X,\beta)$ is the moduli 
stack of stable maps with $m$ parametrizations and no marked points.
This stack admits the action of the torus
$(\C^*)^m$ via its action on $(\P[1])^m$. To carefully give
this torus action, we need some more precise notation. 

\medskip

Fix vector spaces
$W_i \cong \C^2$ for $i=1$,\dots,$m$. On $W_i$, choose coordinates
$x_i,y_i\in W_i^*$, and fix the action of
$\C^*=\C^*_i$ via
\begin{equation*}
  \mu_i\cdot(x_i,y_i) = (x_i,\mu_i y_i).
\end{equation*}
Let $0_i:=(0:1)$ and $\infty_i:=(1:0) \in \P[1]_i$ be the fixed
points of the action of $\C^*_i$ on $\P[1]_i = \P(W_i)$.  Let $\T :=
\prod \C_i^*$ acting diagonally on $\prod \P[1]_i$ and hence on each 
of the graph spaces $\G{m}(X,\beta)$.

\begin{ISC}  The two-parametrized graph space of a point 
\begin{equation*}
  \G{2}(\mathrm{pt},0) = \M{0}(\P[1]\times \P[1],(1,1)) \iso \P[3].
\end{equation*}
A stable map to $\P[1]\times \P[1]$ either embeds $C$ as 
a smooth curve of type $(1,1)$ or as a pair of intersecting rulings.
Thus the stable map space is the linear series.
\end{ISC}

\begin{observation} A stable parametrized map 
$[f] \in \G{m}(X,\beta)$ is a fixed point for the
action of $\T$ described above exactly when:
\begin{itemize}
\item $f$ is constant on each parametrized component

\item Each parametrized component is uniquely parametrized

\item Each node on a parametrized component is at $0$ or~$\infty$.
\end{itemize}
\end{observation}

For our purposes, we will only need to consider the following
``types'' of fixed loci for the action of $\T=\prod \C^*_i$
on the graph space $\G{m+1}(X,\beta)$:

\subsection{Type 1}  
\begin{figure}[ht]
\input{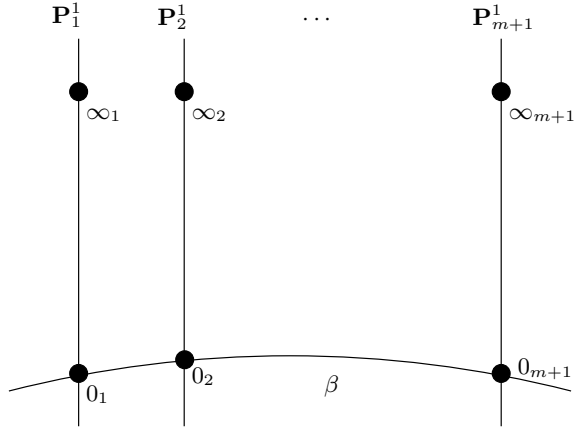}
\caption{Type 1 fixed locus}\label{fig:1}
\end{figure}

A single copy of $\M{m+1}(X,\beta)$ ``embedded
at zeroes''.  (See Figure~\ref{fig:1}.) 

Let $Y = X \times \prod_{i=1}^{m+1} \P[1]_i$, and consider the
gluing morphism:
\begin{equation*}
  \M{m+1}(Y,\beta) \times_Y
  \prod_{i=1}^{m+1}\M{1}(Y,1_i) \to \G{m+1}(X,\beta),
\end{equation*}
where $\beta = (\beta,0^{m+1}) \in \H_2(Y,\Z)$ and likewise
for~$1_i$.
Each $\M{1}(Y,1_i) \iso Y$ and $\M{m+1}(Y,\beta) \cong \M{m+1}(X,\beta) \times
\prod \P[1]_i$, and we obtain a regular embedding
\begin{equation*}
i_{[m],\beta}\colon F_{[m],\beta} := \M{m+1}(X,\beta)
\hra \G{m+1}(X,\beta)
\end{equation*}
by embedding $\M{m+1}(X,\beta) \times
\prod 0_i \hookrightarrow \M{m+1}(X,\beta) \times
\prod \P[1]_i$ and using the gluing morphism above to
further embed in the graph space. 

\subsection{Type 2} 
\begin{figure}[ht]
  \input{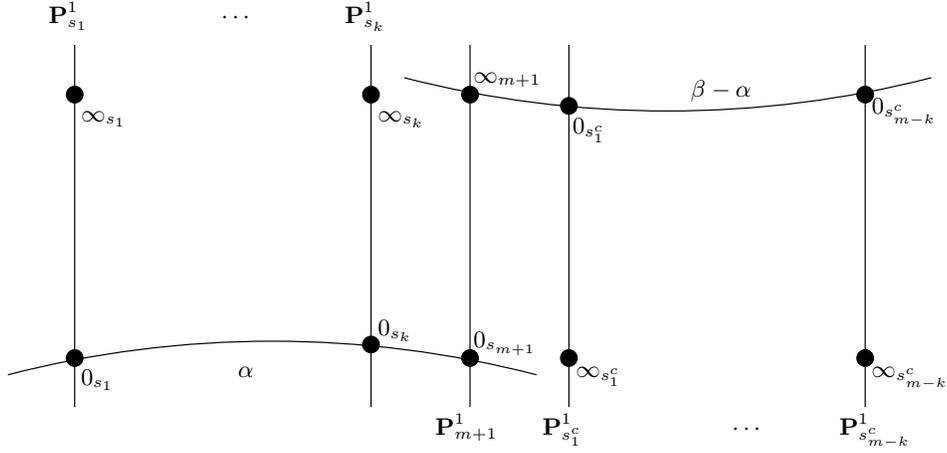}
  \caption{Type 2 fixed locus}\label{fig:2}
\end{figure}

$\M{k+1}(X,\alpha) \times_X
\M{m-k+1}(X,\beta - \alpha)$ ``with $\P[1]_{m+1}$ in the
middle.'' (Copies indexed by subsets $1\in S \subseteq [m]$ with $|S| =
k$ and  $\alpha \in \H_2(X,\Z)$).  (See Figure~\ref{fig:2}.) 

In this case, we consider the composition of gluing maps taking:
\begin{multline*}
  \prod_{s_i\in S}\M{1}(Y,1_i)\times_Y\M{k+1}(Y,\alpha)  
  \times_Y \M{2}(Y,1_{m+1}) \times_Y \\
  \times_Y\M{m-k+1}(Y,\beta - \alpha) \times_Y \prod_{s^c_i\in
  S^c}\M{1} (Y,1_i)
  \to \G{m}(X,\beta),
\end{multline*}
assuming $(S,\alpha) \ne (\{1\},0)$, $([m],\beta)$ or any $(\hat
j,\beta)$. (These appear as other types!)  The product is
isomorphic to
\begin{equation*}
  \prod_{s_i\in S} \P[1]_{s_i} \times 
  \P[1]_{m+1}\times \M{k+1}(X,\alpha) \times_X
  \M{m-k+1}(X,\beta - \alpha) \times \P[1]_{m+1} \times
  \prod_{s^c_i}\P[1]_{s^c_i}
\end{equation*} 
and we identify the embedding: 
\begin{equation*}
  i_{S,\alpha}\colon F_{S,\alpha} = \M{k+1}(X,\alpha) \times_X
  \M{m-k+1}(X,\beta - \alpha) \to \G{m+1}(X,\beta)
\end{equation*}
with $(0_S,0_{m+1})\times \M{k+1}(X,\alpha) \times_X
\M{m-k+1}(X,\beta - \alpha) \times (\infty_{m+1},0_{S^c})$.

\subsection{Type 3} $\M{m}(X,\beta)$ ``with
$\P[1]_{m+1}$ in various places'' (three subtypes). 

\begin{enumerate}
\item $\P[1]_1$ as a tail off of $\P[1]_{m+1}$ (a single copy). Let
\[ 
  i_{\{1\},0}\colon F_{\{1\},0} := \M{m}(X,\beta) \to \G{m+1}(X,\beta)
\]
be the embedding associated to Figure~\ref{fig:3a}. 

\begin{figure}[ht]
  \input{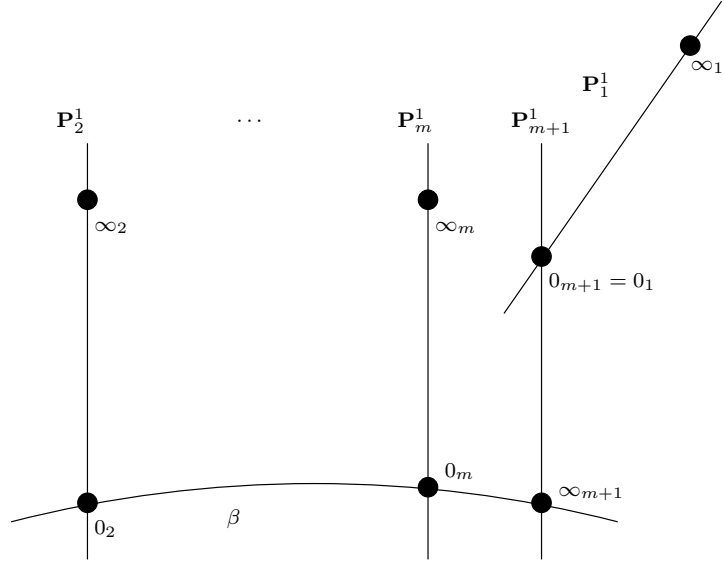}
  \caption{Type 3a fixed locus}\label{fig:3a}
\end{figure}

\item $\P[1]_j$ as a tail off of $\P[1]_{m+1}$ (one
for each $1 < j \le m$).  Let
\[
  i_{\hat j,\beta}\colon F_{\hat j,\beta} := \M{m}(X,\beta) \to
  \G{m+1}(X,\beta)
\]
be the embedding associated to Figure~\ref{fig:3b}. 

\begin{figure}[ht]
  \input{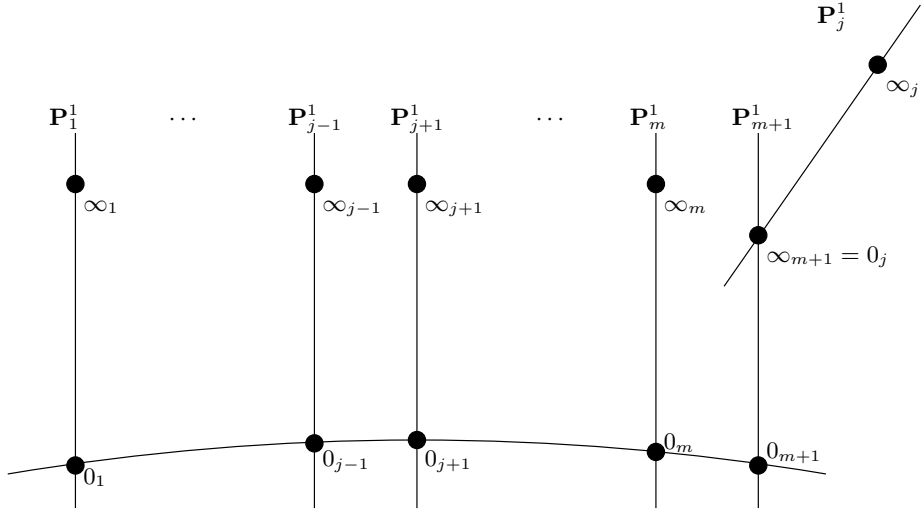}
  \caption{Type 3b fixed locus}\label{fig:3b}
\end{figure}

\item $\P[1]_{m+1}$ as a tail off of $\P[1]_j$
(indexed by $1 < j \le m$).  Let
\[
  i_j\colon F_j := \M{m}(X,\beta) \to \G{m+1}(X,\beta)
\]
be the embedding associated to Figure~\ref{fig:3c}.

\begin{figure}[ht]
  \input{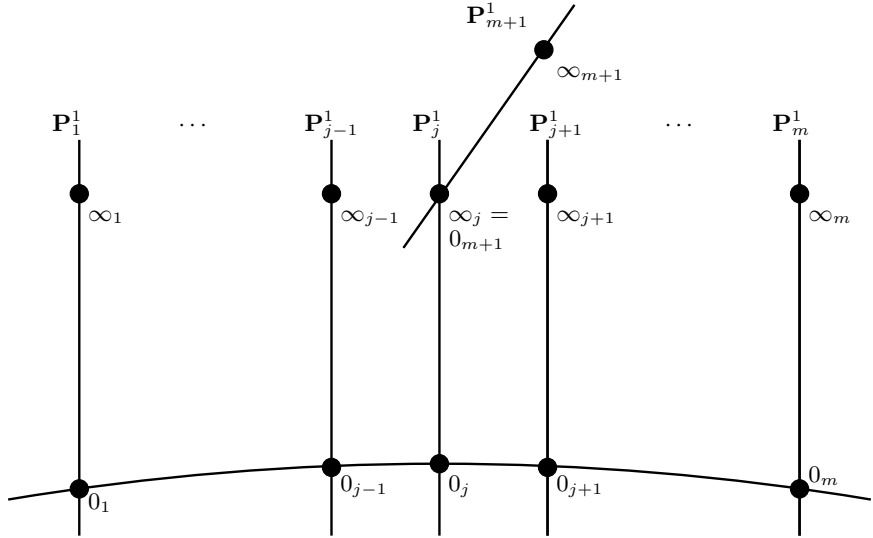}
  \caption{Type 3c fixed locus}\label{fig:3c}
\end{figure}
\end{enumerate}

\begin{lemma}\label{lem:one} There is a $\T$-equivariant birational
morphism: 
\[
  \Phi\colon \G{m+1}(X,\beta) \to \G{m}(X,\beta) \times \P[3]
\]
which (when projected onto the first factor) forgets the last
parametrization and (when projected onto the second factor) forgets
the map to $X$ and all parametrizations except for the first and last.
\end{lemma}

\begin{proof} The existence of $\Phi$ follows from the functoriality
of Kontsevich-Manin spaces. The two projections are just the two
maps:
\[
  \M{0}(X\times \prod_{i=1}^{m+1} \P[1]_i,(\beta,1^{m+1})) 
  \to \M{0}(X\times \prod_{i=1}^{m}\P[1]_i,(\beta,1^{m}))
\] 
and 
\[
  \M{0}(X\times \prod_{i=1}^{m+1} \P[1]e _i,(\beta,1^{m+1})) \to
  \M{0}(\P[1]_1\times \P[1]_{m+1},(1,1))
\]
which clearly commute with the action of~$\T$.

Over the open subset of $\P[3]$ consisting of smooth 
curves, $\Phi$ is an isomorphism. The last parametrization is 
of the same component as the first, and is given by the correspondence
$\P[1]_1 \isoa \P[1]_{m+1}$ induced by the curve
in~$\P[1]_1\times\P[1]_{m+1}$.  
\end{proof}

\begin{lemma}\label{lem:two} Let $(0,0) \in \G{2}(\mathrm{pt},0) = 
\P[3]$ be the  
fixed point corresponding to the  ``coordinate axes.''  Then the
embeddings of Types 1--3 listed above are a complete list of the
fixed loci that are contained in
\[
  \Phi^{-1}\left(F_{[m-1],\beta} \times
  (0,0)\right) \subset \G{m+1}(X,\beta).
\]
Moreover, the induced maps $\Phi|_F\colon F \to F_{[m-1],\beta}\cong
\M{m}(X,\beta)$ are: 

\begin{align*}
  \tag{Type 1} \pi_{m+1}\colon \M{m+1}(X,\beta) & \to
   \M{m}(X,\beta)\\
  \tag{Type 2} \delta_{S,\alpha} \colon \M{k+1}(X,\alpha) \times _X
   \M{m-k+1}(X,\beta - \alpha) & \to \M{m}(X,\beta)\\ 
  \tag{Type 3} \M{m}(X,\beta) & \xra{\id} \M{m}(X,\beta).
\end{align*}
\end{lemma}

\begin{proof} Given a stable map $f\colon C\to X$, 
represent $C$ by a tree with vertices and edges corresponding to
the nodes and components of~$C$, respectively. For an
$f \in \G{m+1}(X,\beta)$ to
map to $F_{[m-1],\beta} \subset \G{m}(X,\beta)$ under the
forgetful map, each $\P[1]_i$ must parametrize a different
curve (edge of the tree) mapping with degree $0$ to $X$, and each
$0_i$ must be a node (vertex) for $i =1$,\dots,$m$. Also, the
shortest path between two such vertices of the tree cannot
contain any such edges, and if one of those edges is removed,
the tree either stays connected or it has two components, one of
which is an edge corresponding to $\P[1]_{m+1}$ mapping with
degree $0$ to~$X$. 
In order for $f$ to map to $(0,0)$ under the
other forgetful map to~$\P[3]$, 
$\P[1]_1$ and $\P[1]_{m+1}$ must represent different
edges, $0_1$ and
$0_{m+1}$ must represent vertices (possibly the same one) and
the shortest path from $0_1$ to $0_{m+1}$ may not contain either
of the two edges. 

The only fixed points under the torus action which satisfy both
conditions are those of
types 1--3. This proves the first part of the lemma.

\medskip

Under the map to $F_{[m-1],\beta}$, the parametrization
of $\P[1]_{m+1}$ is forgotten.  This may result 
in an unparametrized component with $1$ or $2$ nodes, which
is then collapsed.  Moreover, in the $1$ node case, the
resulting marked point must also be forgotten.  This gives the 
second part of the lemma.
\end{proof}

\section{Localization and the main theorem for $\P[n]$.}
When a complex Lie group
$G$ acts on a complex manifold~$X$, there is an
equivariant cohomology ring:
\[
  H^*_G(X,\Q)
\]
which is an algebra over the cohomology ring of the classifying
space~$BG$. If $G = \T = (\C^*)^m$, then
$\H^*(BG,\Q) \cong \Q[t_1,...,t_m]$. The equivariant
cohomology ring for a trivial action of $\T$ is the polynomial 
algebra $\H^*(X,\Q)[t_1,\dots,t_m]$, but in general it is more
complicated. Linearized vector bundles
$E$ on $X$ have equivariant Chern classes $c_d^G(E)$ taking 
values in equivariant cohomology, and equivariant cohomology
pulls back and pushes forward (for proper maps).   

\medskip

The set-up is similar for smooth Deligne-Mumford stacks. In this
setting, there is an equivariant Chow ring (see \cite{EG})
$\CH^*_G(X)$ which always pulls back and pushes forward under
equivariant proper maps. A linearized vector bundle $E$ in this
setting has equivariant Chern classes 
$c_d^G(E) \in \CH^*_G(X)$.

\medskip

A basic result in either setting is the theorem of 
Atiyah-Bott (see \cite{Kr}):

\begin{thm}[Localization]  Each fixed substack
$i\colon F \hookrightarrow X$ of a torus action on a proper smooth
Deligne-Mumford stack is a regularly embedded proper smooth
Deligne-Mumford stack, its normal bundle is  canonically
linearized, its Euler class
$\epsilon_{\T}(F)$ (the top equivariant Chern class of the normal
bundle) is invertible in
\[
  \CH^*(F,\Q)\otimes_{\Q}\Q(t_1,\dots,t_m),
\]
and
any element $c\in \CH_{\T}^*(X)$ is uniquely recovered (modulo torsion) via
the following \emph{localization formula:}
\[
  c = \sum_F i_*\frac{i^*c}{\epsilon_{\T}(F)}.
\]
\end{thm}

Our main interest is in the following simple corollary:

\begin{cor}[Correspondence of residues] Suppose $f\colon X \to X'$
is a $\T$-equivariant map of smooth proper Deligne-Mumford stacks
with $\T$-actions.  If $i'\colon F' \hookrightarrow X'$ is a fixed
substack and $c \in \CH^*_{\T}(X)$, let $f_F\colon F \to F'$ be the
restriction of $f$ to each of the fixed substacks $F\subset
f^{-1}(X')$. Then
\[
  \sum_{F \subset f^{-1}(F')}
  {f_F}_*\frac{i^*c}{\epsilon_{\T}(F)} = \frac{{i'}^*f_*c}
  {\epsilon_{\T}(F')}
\]
\end{cor}

\begin{proof} The two sides of the formula represent the contribution
of $F'$ to localization formulas for $f_*c$ which, by 
uniqueness, must coincide. 
\end{proof}

To prove the main theorem for general $X$
we will need virtual classes.  For now we will prove it
in the case $X = \P[n]$, where the basic idea and most of the
computations are the same, and are not obscured by the presence of
virtual classes.  

\begin{proof}[Proof of the main theorem for~\mbox{$\P[n]$}] Let $c \in
\CH^*_{\T}(\G{m+1}(\P[n],d))$.  Then applying
correspondence of residues to the map $\Phi$ of smooth
Deligne-Mumford stacks (here we use $X = \P[n]$) of Lemma~\ref{lem:one},
and using the enumeration of fixed loci in Lemma~\ref{lem:two}, we get:
\begin{multline}\label{eq:corr}
 {\pi_{m+1}}_*\left(\frac{i_{[m],d}^*c}{\epsilon_{\T}(F_{[m],d})}\right) +  
  \sum_{1\in S}\sum_{e=0}^d{\delta_{S,e}}_*
  \left(\frac{i_{S,d-e}^*c}{\epsilon_{\T}(F_{S,d-e})}\right)
  + \frac{i_{\{1\},0}^*c}
  {\epsilon_{\T}(F_{\{1\},0})} + \\
 \sum_{j=2}^{m}\left(\frac{i_{\hat j,d}^*c}
  {\epsilon_{\T}(F_{\hat j,d})} + \frac{i_{j}^*c}
  {\epsilon_{\T}(F_{j})}\right) = \frac{i_{[m-1],d}^*\Phi_*c}
  {\epsilon_{\T}(F_{[m-1],d})} \cdot \frac 1{t_1t_{m+1}(t_1+t_{m+1})}.
\end{multline}
(The computation 
$\epsilon_{\T}(0,0) = t_1t_{m+1}(t_1+t_{m+1})$ is easily made.)

\medskip

Now, the equivariant Euler classes $\epsilon_{\T}(F)$ appearing in the
denominators depend entirely on the nodes of the domain of a
general representative~$f\in F$.  Essentially, there are two types of
nodes:  those of type~I, where at the point $p_i\in\{0_i,\infty_i\}$,
the $i$th parametrized component meets a component mapping in positive
degree $\alpha$ to $X$, and those of type~II, where at the point
$p_i\in\{0_i,\infty_i\}$, the $i$th parametrized component meets the
point $p_j\in\{0_j,\infty_j\}$ of the $j$th parametrized component.
See Figure~\ref{fig:nodes}.

\begin{figure}[ht]
  \input{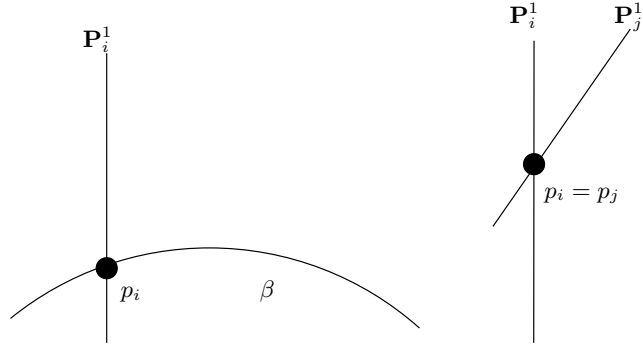}
  \caption{Type I and type II nodes}\label{fig:nodes}
\end{figure}

Any type~I node is a codimension~$2$ condition---one for the node and
one for specifying~$p_i$---while a type~II node is a codimension~$3$
condition---one for the node and two more for specifying $p_i$
and~$p_j$.  Set
\[
  \nu_i:=\begin{cases} 1 & \text{if $p_i=0_i$} \\ 
                      -1 & \text{if $p_i=\infty_i$} \end{cases}.
\]
Then the type~I node contributes the factor 
\[
  \nu_it_i(\nu_it_i-\psi_i)
\]
to $\epsilon_{\T}(F)$, while the type~II node contributes
\[
  \nu_it_i\nu_jt_j(\nu_it_i+\nu_jt_j).
\]
(For this type of computation, see \cite{BCPP}.)

Thus, if we let 
\[
  t := t_{m+1},
\] 
the following computations are valid on any $\G{m+1}(X,\beta)$:
\addtocounter{equation}{1}
\begin{align*}
 \tag{\arabic{equation}.1}
  \epsilon_{\T}(F_{[m],\beta}) & = 
  t(t-\psi_{m+1})\prod_{i\in [m]} t_i(t_i - \psi_i).\\ 
 \tag{\arabic{equation}.2}
  \epsilon_{\T}(F_{S,\alpha}) & = t(t-\psi_{k+1})(-t)(-t - \psi_1')
  \prod_S t_{s_i}(t_{s_i} - \psi_i)
  \prod_{S^c} t_{s^c_i}(t_{s^c_i} - \psi_{i+1}') \\
\intertext{where $\psi_i,\psi_i'$ are the cotangent classes on
$\M{k+1}(X,\alpha)$ and $\M{m-k+1}(X,\beta - \alpha)$.}
 \tag{\arabic{equation}.3a}
  \epsilon_{\T}(F_{\{1\},0}) & = t_1t(t_1+t)(-t)(-t -
  \psi_1)\prod_{i=2}^m t_{i}(t_i - \psi_i).\\ 
 \tag{\arabic{equation}.3b}
  \epsilon_{\T}(F_{\hat j,\beta}) & = (-t)t_j(-t+t_j)\prod_{S=\hat
  j}t_{s_i}(t_{s_i} - \psi_i)t(t-\psi_m).\\
 \tag{\arabic{equation}.3c}
  \epsilon_{\T}(F_j) & = (-t_j)t(-t_j+t)\prod_{S=\hat
  j}t_{s_i}(t_{s_i} - \psi_i)t_j(t_j - \psi_m).
\end{align*}

Finally, let 
$\P[n] = \P(V)$ and consider $H_{\T}$, 
the equivariant hyperplane class on the linear space
$\P(\Hom(\Sym^d(W_{m+1}),V))$.  There is an 
equivariant morphism 
\[
   v\colon \G{m+1}(\P[n],d) \to \G{1}(\P[n],d) \to
  \P(\Hom(\Sym^d(W_{m+1}),V)) 
\]
which is a composition of the forgetful map remembering only 
$\P[n]$ and the last parametrization, and the ``map to the
linear sigma model.''   (The geometry of this second (birational) map
was used in \cite{B1} to give a proof of the mirror theorem.)  This
$H_{\T}$ pulls back to the fixed loci as follows:

\addtocounter{equation}{1}
\begin{align*}
 \tag{\arabic{equation}.1}
  i_{[m],d}^*v^*H_{\T} & = e_{m+1}^*H\\ 
 \tag{\arabic{equation}.2}
  i_{S,d-e}^*v^*H_{\T} &= e_{k+1}^*H - et\\
 \tag{\arabic{equation}.3a}
  i_{\hat j,d}^*v^*H_{\T} & = e_j^*H\\
 \tag{\arabic{equation}.3b}
  i_{\{1\},0}^*v^*H_{\T} & = e_1^*H - dt\\
 \tag{\arabic{equation}.3c}
  i_{j}^*v^*H_{\T} & = e_j^*H.
\end{align*}

To see this, note that under the morphism~$v$, the fixed loci map
to various copies of $\P[n]$ sitting as the fixed loci in 
$\P(\Hom(\Sym^d(W_{m+1}),V))$.  More specifically, set
\begin{multline*}
  (\P[n])_e := \{x_{m+1}^{d-e}y_{m+1}^e\} \times \P(V) \subset  
   \P\big(\Sym^d(W_{m+1}^*)\big) \times \P(V) \xra{\mathrm{Segre}}\\
   \P(\Hom(\Sym^d(W_{m+1}),V)) 
\end{multline*}
which are all fixed under the $\T$ action.  One easily computes that
$H_{\T}$ restricts to $\H^*_{\T}((\P[n])_e,\Q)$ as~$H-et$.   Finally,
fixed loci of types 1, 3a, and 3c map under $v$ to $(\P[n])_0$, the
fixed loci of type 3b map to $(\P[n])_d$, and the loci of type 2 map to
$(\P[n])_e$ for appropriate $1\le e\le d-1$.

Substitute in the summands of \eqref{eq:corr} for the equivariant Euler
classes and for the choice $c = v^*H_{\T}^b$ and
push forward under the total evaluation map~$\ev$.  Then by the
projection formula and the computations above we obtain the
following:
\addtocounter{equation}{1}
\begin{align*}
 \tag{\arabic{equation}.1}
  \ev_*{\pi_{m+1}}_*\frac{i_{[m],d}^*c}{\epsilon_{\T}(F_{[m],d})} & = 
  J^{\P[n]}_d(t_1,\dots,t_m,t)\otimes_{H^b} J^{\P[n]}_0(-t) \\
 \tag{\arabic{equation}.2}
  \ev_*{\delta_{S,\alpha}}_*\frac{i_{S,\alpha}^*c}{\epsilon_{\T}(F_{S,\alpha})}
  & =  J^{\P[n]}_{d-e}(\vec t_S,t)\otimes_{(H - et)^b}
  J^{\P[n]}_e(-t,\vec t_{S^c})\\
 \tag{\arabic{equation}.3a}
  \ev_*\frac {i_{\{1\},0}^*c}{\epsilon_{\T}(F_{\{1\},0})} & =
   J^{\P[n]}_0(t_1,t) \otimes_{(H-dt)^b} J^{\P[n]}_d(-t,\vec
   t_{\hat1})\\
 \tag{\arabic{equation}.3b}
  \ev_*\frac {i_{\hat j,d}^*c}{\epsilon_{\T}(F_{\hat j,d})} & = 
   J^{\P[n]}_d(\vec t_{\hat j},t)\otimes_{H^b}
   J^{\P[n]}_0(-t,t_j) \\
 \tag{\arabic{equation}.3c}
  \ev_*\frac {i_j^*c}{\epsilon_{\T}(F_j)} & = J^{\P[n]}_d(\vec t_{\hat
  j},t_j)\otimes_{H^b}J^{\P[n]}_0(-t_j,t)
\end{align*}

We get the theorem by multiplying both sides of \eqref{eq:corr} by
$t_1t(t_1+t)$, collecting types 1, 2, 3a, and 3b under the 
double sum, and  noting that
\[
  \ev_*\frac{i_{[m-1],d}^*\Phi_*c}{\epsilon_{\T}(F_{[m-1],d})} \in
  \H^*((\P[n])^m,\Q)[t_1,t_1^{-1},\dots,t_m,t_m^{-1},t]
\]
because $i_{[m-1],d}^*\Phi_*c$ is {\bf polynomial} in
$t_1$,\dots,$t_m,t$ and the inverse to $\epsilon_{\T}(F_{[m-1],d})$
belongs to $\H^*(\M{m}(\P[n],d),\Q)[t_1^{-1},\dots,t_m^{-1}]$.
\end{proof}

\section{Virtual classes}
In order to prove our main theorem in the general case, we will
need to establish a simple property of equivariant virtual classes.
We begin by briefly recalling the
construction of the virtual class on Kontsevich-Manin spaces, following
Behrend and Fantechi~\cite{Beh,BF} (but see also Li-Tian~\cite{LT}) and
of the equivariant virtual class on graph 
spaces, following Graber-Pandharipande~\cite{GP}.

\medskip

Fix a complex projective manifold $X$ and an embedding $X \hra
\P[n]$.  For each $\beta\in \H_2(X,\Z)$, 
let $d$ be the degree of the image of $\beta$ in~$\P[n]$. Then
there  is a commuting diagram of stacks:
\begin{equation*}
\begin{matrix}
\M{m}(X,\beta) & & \xra{i} & & \M{m}({\bf P}^n,d) \\ 
& \searrow & & \rho\swarrow \\ 
& & \goth{M}_{0,m}
\end{matrix}
\end{equation*}
where $\goth{M}_{0,m}$ is the Artin (not Deligne-Mumford!) stack of
prestable $m$-marked curves. The map $i$ is a closed
embedding (let $\sh{I}$ be the associated ideal sheaf) and $\rho$
is smooth. It follows that the relative intrinsic
normal cone of Behrend-Fantechi is the  cone stack associated to the
following map of sheaves on $\M{m}(X,\beta)$:
\begin{equation*}
  \sh I/\sh I^2 \to i^*\Omega_{\M{m}(\P[n],d)/ \goth M_{0,m}}
\end{equation*}

This relative normal cone, which we denote by
$\goth{C}_{\M{m}(X,\beta)/\goth{M}_{0,n}}$, embeds in the 
smooth $h^1/h^0$ cone stack $V_{\M{m}(X,\beta)/\goth{M}_{0,m}}$
associated to the object
\[
(R^1\pi_*e^*TX)^\vee
\]
of the derived category of coherent sheaves on $\M{m}(X,\beta)$.
(The dual is the Verdier dual and $\pi\colon{\mathcal C} \to
\M{m}(X,\beta)$ and $e\colon {\mathcal C} \to X$ come from the 
universal stable map). 
The virtual class $\vir{\M{m}(X,\beta)}$ is then obtained by 
pulling back the class of $\goth{C}_{\M{m}(X,\beta)/\goth{M}_{0,m}}$
via the zero section of $V_{\M{m}(X,\beta)/{\goth M}_{0,m}}$.

\medskip

Similarly, for graph spaces, there is a diagram of $\T$-invariant
morphisms
\[
\begin{matrix} \G{m}(X,\beta) & & \xra{i}
& & 
\G{m}(\P[n],d) \\ & \searrow & & \rho\swarrow \\ 
& & \goth{G}_{0,m}
\end{matrix},
\]
where $\goth{G}_{0,m}$ is the stack
of prestable zero-pointed maps to $(\P[1])^m$ of 
multi-degree $(1,\dots,1)$.
The (equivariant) intrinsic relative normal cone
$\goth{C}_{\G{m}(X,\beta)/\goth{G}_{0,m}}$ and $h^1/h^0$  
cone $V_{\G{m}(X,\beta)/\goth{G}_{0,m}}$ are defined exactly
as before, and the (equivariant) virtual class 
$\vir{\G{m}(X,\beta)}_\T$ may also be defined as before, using
equivariant Chow groups. This definition is simpler than the
definition in \cite{GP}, but is equivalent. The simplification in 
our case comes from the existence of the $\T$-invariant embedding $i$
into a relatively smooth graph space.
The significance of the equivariant virtual class is in the following
``virtual'' version of the localization theorem:

\medskip

\begin{thm}[Graber-Pandharipande] In the equivariant Chow group of 
the graph space $\G{m}(X,\beta)$ the virtual class satisfies
\[
  \vir{\G{m}(X,\beta)}_{\T} = \sum_F
  i_*\frac{i^*\vir{\G{m}(X,\beta)}_{\T}}{\epsilon_{\T}(F)} 
\]
where $i\colon F \hra\G{m}(X,\beta)$ are the (regular)
embeddings of the fixed substacks.
\end{thm}

In order to use this theorem, we need the following:
\begin{lemma}\label{lemma:X} \ 

\begin{enumerate}
\item\label{num:Xi} The forgetful map 
\[
  \phi\colon \G{m}(X,\beta) \to \M{0}(X,\beta)
\]
is flat and equivariant for the trivial action of $\T$
on~$\M{0}(X,\beta)$.

\item\label{num:Xii} The equivariant virtual class satisfies
\[
 \vir{\G{m}(X,\beta)}_{\T} = \phi^*\vir{\M{0}(X,\beta)},
\]
where $\vir{\M{0}(X,\beta)}$ is the ordinary virtual
class, regarded as an equivariant class for the trivial action
of~$\T$. In particular, each 
$i^*\vir{\G{0}(X,\beta)}_{\T} = \vir{F}$ in the theorem above, where 
$\vir{F}$ is the ``ordinary'' virtual class on~$F$, thought of
as a fiber product of Kontsevich-Manin spaces.
\end{enumerate}
\end{lemma}

\begin{proof} It suffices by induction to prove the lemma for the case
$m=1$. In that case, we will consider a (non-commuting!) diagram of
stacks:
\[
\begin{CD} \M{3}(X,\beta) @>{g}>> \G{1}(X,\beta) @>{\phi}>> \M{0}(X,\beta)
\\
@VVV @VVV @.\\
\goth{M}_{0,3} @>{\goth{g}}>> \goth{G}_{0,1} @. {}
\end{CD}
\]
where the horizontal maps are the ``cross-ratio'' maps defined as
follows. The universal curve ${\mathcal C} \cong \M{4}(X,\beta)$
over $\M{3}(X,\beta)$
maps to $\M{4}\cong\P[1]$ via the forgetful map. Together with 
the evaluation map to $X$, this defines~$g$. If $f\colon C\to X$ is 
a stable map with $3$ marked points $p$, $q$, $r\in C$, then 
the map $C \to \P[1]$ defined by $g(f)$ may be taken to 
be the unique map with the property that $f(p) = 0, f(q) = 1$ and
$f(r) = \infty$. This is the cross-ratio if $p$, $q$, $r$ belong to the 
same component of~$C$, but is well-defined even if they lie on
different components.

\medskip

For $\goth{g}$, we apply the prestabilization map ${\mathcal C}
\to  \goth{M}_{0,4}$ (see \cite{Beh}) to the universal curve over
$\goth{M}_{0,3}$ followed by the stabilization map
$\goth{M}_{0,4} \to \M{4} \cong \P[1]$. This map has
the same pointwise description as~$g$.

\medskip               

The diagram doesn't commute because $g$  
stabilizes unstable maps to $X\times
\P[1]$, while $\goth{g}$ does not. On the
other hand, there is a ``good'' open substack
\[
U := \{ f:C \to \P[1] \mid \text{$f$ is an isomorphism over $0$,
$1$, $\infty$}\} \subset \goth{G}_{0,1}
\]
with the following properties:

\begin{itemize}
\item $g$ and ${\goth g}$ are both isomorphisms over $U$.

\item The diagram above is Cartesian when restricted to~$U$.

\item Translates of $U$ by elements $m\in \mathrm{PGL}(2,\C)$ 
cover~$\goth{G}_{0,1}$,
\end{itemize}
\medskip

If $f\in U$, then $p$, $q$, and $r$ are the preimages of $0$, $1$, and
$\infty$, so ${\goth g}$ is invertible at~$f$. If $f\in \M{3}(X,\beta)$ lies 
over $U$, then $p$, $q$, and $r$ all belong to same 
component $C_0\subset C$ of the curve associated to $f$, and $g(f)$
imposes the unique parametrization on $C_0$ taking $p$, $q$, and $r$ to
$0$, $1$, and~$\infty$. Clearly, then, $g$ and ${\goth g}$ are isomorphisms
over $U$ and the diagram is Cartesian over~$U$. Since every prestable
map $f\colon C \to \P[1]$ of degree one is {\it generically} an
isomorphism over~$\P[1]$, it follows that the translates of $U$ cover
$\goth{G}_{0,1}$. 

\medskip

We finish the proof now by comparing $\G{1}(X,\beta)$ with
$\M{3}(X,\beta)$. Suppose $f\in \G{1}(X,\beta)$ lies over~$U$.
Then $g$ is an isomorphism at $f$, so since $\phi \circ g$ is flat
everywhere (it is a composition of the flat forgetful maps), it
follows that $\phi$ is flat at~$f$. But an arbitrary $f\in
\G{1}(X,\beta)$ lies over some translate $mU$, over which the
composition of $g$ with translation by $m$ is an isomorphism, and
we similarly conclude that $\phi$ is flat at an arbitrary $f$. This 
gives us~(\ref{num:Xi}).

\medskip

Thus $\phi$ is flat, and we may use the flat pull-back to define 
$\phi^*\vir{\M{0}(X,\beta)}$. Behrend showed that 
the relative intrinsic normal cone 
${\goth C}_{\M{0}(X,\beta)/\goth{M}_{0,0}}$ pulls back under $
\phi\circ g$ to ${\goth C}_{\M{3}(X,\beta)/\goth{M}_{0,3}}$ and the
same trick we employed in the previous paragraph shows that it pulls
back under
$\phi$ to ${\goth C}_{\G{1}(X,\beta)/\goth{G}_{0,1}}$. The flatness
of 
$\phi$ also tells us that $R^1\pi_*e^*TX$ pulls back to the 
corresponding element of the derived category of sheaves on
$\G{1}(X,\beta)$, and we get~(\ref{num:Xii}). The last sentence of
(\ref{num:Xii}) is  
a consequence of Behrend's work, since the induced maps
$F \rightarrow \M{0}(X,\beta)$ are always gluing maps of
Kontsevich-Manin spaces.
\end{proof}
 
\section{The main theorem and reconstruction.} 
We now return to Theorem~\ref{thm:main1} and its generalizations and
consequences. 

\begin{proof}[Proof of the main theorem (rank one case):]  We may assume
that $H$ is very ample. Indeed, suppose the
polynomiality condition holds for the expression
\[
  \sum_{1\in S \subseteq [m]}\sum_{e=0}^d
  J^X_{d-e}(\vec t_S,t)\otimes_{(lH-et)^b}
  J^X_e(-t,\vec t_{S^c}) + \sum_{j=2}^m
  J^X_d(\vec t_{\hat j},t_j)\otimes_{(lH)^b}
 J^X_0(-t_j,t)
\]
for some $l > 0$. Only the $e$'s divisible by $l$ will
produce non-zero terms, because the degree of every curve
(measured against $lH$) is a multiple of $l$.  But 
replacing $lH - et$ by
$H - \frac elt$ in the twisted tensor products simply multiplies
the expression by~$l^{-b}$. If we now replace the subscript of 
each $J$ by
the  degree of the curve against $H$ (instead of against $lH$) we
get the desired result for~$H$.  

\medskip

The embedding $X \subset \P[n]$ defined by $H$ allows us to
define a morphism
\[
  v\colon \G{m+1}(X,d) \to \G{1}(X,d)\hra \G{1}(\P[n],d)
  \to \P(\Hom(\Sym^d(W_{m+1})\otimes V))
\]
and an equivariant Chern class $v^*(H_{\T}^b)$ as in the 
$\P[n]$ case. Applying Lemma~\ref{lemma:X}~(\ref{num:Xi}) to the map
$\G{m+1}(X,\beta) \to \G{m}(X,\beta)$, we see that $\Phi$ is a
local complete intersection (l.c.i.) morphism, since it
factors through the graph followed by a flat morphism: 
\[
  \begin{matrix} & & \G{m+1}(X,d) \times \P[3] \\
  & \nearrow & \downarrow \\ 
  \G{m+1}(X,d) & \stackrel{\Phi}
   \to & \G{m}(X,d) \times \P[3] 
\end{matrix}.
\]
Then by Lemma~\ref{lemma:X}~(\ref{num:Xii}),
\[
\Phi^*\left(\vir{\G{m}(X,d)}_{\T} \times [\P[3]]\right) = 
\vir{\G{m+1}(X,d)}_{\T}.
\]
It follows by the projection formula that the
correspondence of residues holds for  
$c \cap \vir{\G{m+1}(X,d)}_{\T}$ (and any equivariant
Chern class~$c$) with each $i^*c$
replaced by
$i^*c \cap \vir{F}$, and
$i_{[m-1],d}^*\Phi_*c$ replaced by $i_{[m-1],d}^*\Phi_*c \cap
\vir{F_{[m-1],d}}$ (again, using Lemma~\ref{lemma:X}).
The proof of the $\P[n]$ case now carries over to 
prove the general rank one case.
\end{proof}

Next we turn to the theorem for arbitrary~$\H^2(X,\Q)$. It
seems best to do this, not for $J$-functions defined intrinsically
on~$X$, but for $J$-functions defined in terms of a choice of
(generalized) polarization on~$X$.  (See also \cite{B2}.)

\begin{defn} \ 

\begin{enumerate}
\item A divisor $H$ on $X$ is {\sl eventually free} if some positive
multiple $lH$ defines a morphism $X \to \P[n]$.

\item\label{num:ample} A collection $H_1$,...,$H_k$ (written $H$ for short) of
eventually free divisors is {\sl ample} if positive $\Z$-linear
combinations $l_1H_1 +\dots + l_kH_k$ are ample.

\item The $J$-functions associated to an $H$ as in (\ref{num:ample}) are
\begin{align*}
 J^{X,H}_d(t_1,\dots,t_m) & =
  J^{X,H_1,\dots,H_k}_{(d_1,\dots,d_k)}(t_1,\dots,t_m)\\
  & := \sum_{d(\beta) = (d_1,\dots,d_k)} J^X_\beta(t_1,\dots,t_m),
\end{align*}
where $d(\beta)$ is the multi-degree $(\deg_{H_1}(\beta),\dots,
\deg_{H_k}(\beta))$.
\end{enumerate}
\end{defn}

\medskip

\begin{thm}[The main theorem---general case] If $X$ is a complex
projective manifold and $H = (H_1,\dots,H_k)$ is an ample collection
of eventually free divisors, then
\begin{multline*}
  t(t_1+t)\Bigg(\sum_{1\in S \subseteq [m]}\sum_{e \preceq d}
   J^{X,H}_{d-e}(\vec t_S,t)\otimes_{\prod (H_i-e_it)^{b_i}}
   J^{X,H}_e(-t,\vec t_{S^c}) + \\
  \sum_{j=2}^m J^{X,H}_d(\vec t_{\hat
   j},t_j)\otimes_{\prod H_i^{b_i}} J^{X,H}_0(-t_j,t)\Bigg) \in \H_*(X^m,\Q)[t_1,t_1^{-1},\dots,t_m,t_m^{-1},t] 
\end{multline*}
In this case, we sum over $0 \preceq e = (e_1,\dots,e_k) \preceq d$, meaning
that $0 \le e_i \le d_i$. 
\end{thm}

\begin{proof}
As in the proof of the rank one version, we may assume that 
$H_1$,\dots,$H_k$ are not just eventually free, but free, by replacing
them with positive multiples (which can be taken to be the same
multiple). The $H_i$ define a morphism
\[
  v\colon \coprod_{d(\beta) = d}\G{m+1}(X,\beta) \to
  \coprod_{d(\beta) = d} \G{1}(X,\beta) \to  
  \prod_{i=1}^k \P(\Hom(\Sym^{d_i}W_{m+1},V_i))
\]
and the theorem results from applying the correspondence of 
residues to the class $v^*\prod_{i=1}^k{H_i}_{\T}^{b_i}$, where the
${H_i}_{\T}$ are the equivariant hyperplane classes pulled back
from~$\P(\Hom(\Sym^{d_i}W_{m+1},V_i))$. 
\end{proof}

Finally, we have the
\begin{thm}[Reconstruction]\label{thm:recon}  Let $R_H \subset \H^*(X,\Q)$
be the subring generated as a $\Q$-algebra by $1$ and an ample
collection $H_1$,\dots,$H_k$ of eventually free divisors. If the
orthogonal complement to $R_H$ annihilates each of the
one-variable $J$-functions 
$J^{X,H}_d(t)$, then the Gromov-Witten invariants of the form
\[
  \sum_{d(\beta) = d}
  \left<\gamma_1\psi^{a_1},\dots,\gamma_m\psi^{a_m}\right>^X_\beta
\]
for $\gamma_i \in R_H$ are completely determined by the one-point
invariants, the intersection matrix on~$R_H$, and the canonical
class~$K_X$. 
\end{thm}

\begin{proof} The only term in the main theorem
involving a $J$-function of $m+1$ variables and curves of
(multi) degree $d$ is
\[
  J^{X,H}_d(t_1,\dots,t_m,t)\otimes_{\prod H_i^{b_i}}J^{X,H}_0(-t) = 
   \pi^X_*\left(\left(\cup \pi_X^*\prod H_i^{b_i}\right) \cap
   J^{X,H}_d(t_1,\dots,t_m,t)\right). 
\] 
The product $t(t+t_1)J^{X,H}_d(t_1,\dots,t_m,t)$
is a polynomial in~$t^{-1}$, expanding as
\[
  t(t+t_1)J^{X,H}_d(t_1,\dots,t_m,t) = 
  (t+t_1) \sum_{a=1}^N t^{-a}\sum_{d(\beta) = d} 
  \ev_*\frac{\psi_{m+1}^{a-1} \cap \vir{\M{m+1}(X,\beta)}}
  {\prod_{i=1}^m t_i(t_i - \psi_i)},
\]
for some $N$ depending on~$K_X$. 
It follows by downward induction on the power of $t^{-1}$ 
and the main theorem that every term in the expansion of
$\pi^X_*(J^{X,H}_d(t_1,\dots,t_m,t) \cup \pi_X^*\prod H_i^{b_i})$ 
in $t^{-1}$ is determined inductively by $J$-functions involving fewer
variables and/or lower degrees.  Note that by stopping the induction
at the $t^{-1}$ term, we  determine the constant term, about which the  
main theorem tells us nothing.

\medskip

This argument only proves the reconstruction theorem when all 
cohomology is generated by the $H_i$ since it (inductively)
requires knowledge of the classes 
$\pi^X_*((\pi_X^*\gamma)\cap J^{X,H}_d(t_1,\dots,t_m,t))$ where
$\gamma$ is an \emph{arbitrary} cohomology class. This argument does,
however, capture the main idea of the proof.

\medskip

We now prove the following by induction on $(m+1,d)$:

\begin{claim}\label{cl:1} \  
\begin{enumerate}
\item\label{num:indi} For all $\gamma_1,\dots,\gamma_m \in R_H$ and
$\alpha \in R_H^{\perp}$,
\[
  \sum_{d(\beta) = d} \left<\gamma_1\psi^{a_1}, \dots,
  \gamma_m\psi^{a_m}, \alpha\psi^a\right>^X_\beta = 0.
\]

\item\label{num:indii}  For all $\gamma_1,\dots,\gamma_m \in R_H$,
the invariants
\[
  \sum_{d(\beta) = d} \left<\gamma_1\psi^{a_1}, \dots,
  \gamma_m\psi^{a_m}, \gamma_{m+1}\psi^{a_{m+1}}\right>^X_\beta
\]
are determined by the one-point invariants and the intersection matrix
on~$R_H$. 
\end{enumerate}
\end{claim}
In terms of $J$-functions (using the symmetry), this claim is
equivalent to
\begin{claim}\label{cl:2} \ 
\begin{enumerate}
\item\label{num:ind2i} If $\gamma_1,\dots,\gamma_m \in R_H$ and
$\alpha \in R_H^{\perp}$ then
\[
  \deg \left(\left(\pi_1^*\alpha \cup \pi_2^*\gamma_1 \cup \dots \cup
   \pi_{m+1}^*\gamma_m\right)\cap J^{X,H}_d(t_1,\dots,t_m,t)\right)=0. 
\]
\item\label{num:ind2ii} For all $\gamma_1,\dots,\gamma_m \in R_H$,
\[
  \deg\left(\left( \pi_1^*\gamma_1 \cup \pi_2^*\gamma_2 \cup \dots \cup
   \pi_{m+1}^*\gamma_{m+1}\right)\cap J^{X,H}_d(t_1,\dots,t_m,t) \right)
\]
is determined by one-point invariants and the intersection matrix on~$R_H$.
\end{enumerate}
\end{claim}

To start our induction, note that the claim holds for $m=0$ by
assumption.  Also, the claim holds for $d=0$:
\[
  \deg ((\pi_1^*\alpha \cup \pi_2^*\gamma)\cap J^{X,H}_0(t_1,t_2))
   =\frac 1{t_1t_2(t_1+t_2)} \int_{X} \alpha \cup \gamma = 0
\]
by orthogonality, and
\[
  \deg ((\pi_1^*\gamma_1 \cup \pi_2^*\gamma_2)\cap J_0^{X,H}(t_1,t_2))
    = \frac 1{t_1t_2(t_1+t_2)}\int_X
   \gamma_1 \cup \gamma_2
\]
and hence is determined by the intersection matrix on~$R_H$. 

\medskip

Using the argument at the beginning of this proof, the vanishing
in Claim~\ref{cl:2}(\ref{num:ind2i})  will follow by induction (on the
power of  $t^{-1}$), once we establish vanishing for all
expressions of the form
\begin{align*}
  I_a & :=\deg\left(\left(\pi_1^*\alpha \cup
   \pi_2^*\gamma_1 \cup \dots \cup
   \pi_{m}^*\gamma_{m-1}\right)\cap\left(J^{X,H}_{d-e}(\vec t_S,t)
   \otimes_{\prod (H_i-e_it)^{b_i}} J^{X,H}_e(-t,\vec t_{S^c})\right)\right)\\
\intertext{and}
  I_b & := \deg\left( \left( \pi_1^*\alpha \cup \pi_2^*\gamma_1
   \cup \dots \cup \pi_{m}^*\gamma_{m-1}\right)
   \cap\left(J^{X,H}_d(\vec t_{\hat j},t_j) \otimes_{\prod 
   H_i^{b_i}} J^{X,H}_0(-t_j,t)\right) \right).
\end{align*}
But these expressions may be rewritten:
\[
  I_b = \deg \left(\left(\pi_1^*\alpha
   \cup \pi_2^*\gamma_1 \cup \dots \cup \pi_{m}^*(\gamma_{m-1} \cup
   \prod H_i^{b_i})\right)\cap J^{X,H}_d(\vec t_{\hat j},t_j)\right).
\]
To rewrite $I_a$, choose an orthogonal basis $\lambda_j,\alpha_l\in
\H^*(X,\Q)$ such that $\lambda_j \in R_H$ with intersection matrix
$g_{jj'}$ and $\alpha_l \in R_H^\perp$ with 
intersection matrix $h_{ll'}$. Then
\begin{multline*}
  I_a = \sum_{j,j'}\deg\left(\left(
   \pi_1^*\alpha \cup \dots \cup \pi_{k+1}^*(\prod (H_i - e_it)^{b_i}
   \cup \lambda_j)\right)\cap J^{X,H}_{d-e}(\vec t_S,t) \right)
   g^{jj'}\\ 
  \deg\left(\left( \pi_1^*(\lambda_{j'})\cup \pi_2^*\gamma_k
   \cup \dots \cup 
   \pi_{m-k+1}^*\gamma_{m-1}\right) \cap J^{X,H}_e(-t,\vec t_{S^c})\right)\\
  + \sum_{l,l'} \deg\left(\left( 
    \pi_1^*\alpha \cup \dots \cup \pi_{k+1}^*(\prod (H_i - e_it)^{b_i}
    \cup \alpha_l)\right)\cap J^{X,H}_{d-e}(\vec t_S,t)\right)h^{ll'} \\
  \deg\left(\left(   \pi_1^*(\alpha_{l'})\cup \pi_2^*\gamma_k
   \cup \dots \cup 
   \pi_{m-k+1}^*\gamma_{m-1}\right) \cap J^{X,H}_e(-t,\vec t_{S^c})\right).
\end{multline*}

Now suppose Claim~\ref{cl:2}(\ref{num:ind2i}) holds for all $(n+1,e)$
such that either $n < m$ or $n = m$ and $e \prec d$. Then $I_b=0$
(taking $n = m-1$), and $I_a=0$ since the first factors in the first
double sum and the second factors in the second sum vanish. This proves
Claim~\ref{cl:1}(\ref{num:indi}) by induction.  Similarly, 
assuming Claim~\ref{cl:1}(\ref{num:indi}), we prove
\ref{cl:1}(\ref{num:indii}) by induction,  noting that in this case,
the second double sum in $I_a$ (but not the first) vanishes. The first
double sum and the $I_b$ terms are explicitly determined by the
intersection matrix $g_{jj'}$ and Gromov-Witten invariants for lower~$(n+1,e)$.
\end{proof}

\appendix

\section{Small quantum product for complete intersections} 
We may turn Formula~\ref{form:one} into an algorithm for 
producing structure constants for the small quantum product
on Fano complete intersections in~$\P[n]$.

Given the type, $(l_1,\dots,l_m)$ of the Fano
complete intersection~$S\subset\P[n]$, set
\begin{align*}
f & := n+1-l_1 - \dots -l_m\\
\intertext{the Fano index of $S$, and}
d_{\max} & := \left\lfloor\frac{n-m+1}f\right\rfloor,
\end{align*}
the maximal degree $d$ for which nonzero ``unmixed'' $2$-point
invariants $\left<H^a,H^b\right>^X_d$ may occur (by a dimension count).

\medskip

For $d = 1$,\dots,$d_{\max}$, let $v(d)$ be the vector of one-point
invariants, i.e., $v(d)$ is defined by
\[  
e_*\left(\frac{\vir{\M{1}(X,d)}}{t(t-\psi)}\right) = 
v(d)_0t^{-f} + v(d)_1Ht^{-f-1} + \dots +
v(d)_{n-m}H^{n-m}t^{-f-n+m}.\]
(These are computed by Givental's formulas, Theorem~\ref{thm:giv}.)

\medskip

We define shift matrices of size $(n-m+1)\times (n-m+1)$:
\[
S(d) :=\begin{pmatrix} 
d & 0 & \dots & 0 & 0 \\
1 & d & \dots & 0 & 0 \\ 
& & \vdots \\ 
0 & 0 & \dots & d & 0 \\
0 & 0 & \dots & 1 & d
\end{pmatrix}.
\]
Applying $S(d)$ to a vector corresponds to multiplication by~$(H+dt)$.

\medskip

We define the matrices of mixed invariants, also of size
$(n-m+1)\times (n-m+1)$, indexed from $0$ to $n-m$:
\[
  M(d)_{n-m-a,b} := \frac {(-1)^{c-1}}{\prod l_i}
   \left<H^a,H^b\psi^c\right>^X_d, 
\]
where $c := df +n-m-a-b$.  This is the matrix associated to the
operator
\[
H^b \mapsto {e_1}_*\left(\frac{e_2^*H^b \cap
\vir{\M{2}(X,d)}}{-t-\psi}\right).
\]
It is important to note that $M(d)_{n-m-a,b} = 0$ when $c <0$.

\medskip

In terms of these data structures, our formula becomes a recursive
formula for the $b$th column of $M(d)$ in terms of the lower
$M(e)$'s:
\[
  M(d)_{*,b} = - S(d)^bv(d) - \sum_{e=1}^{d-1}M(d-e)S(e)^bv(e),
\]
except that we must set $M(d)_{n-m-a,b} = 0$ whenever $c < 0$. This
amounts to truncating $M(d)$ at the upper right corner.

\medskip

Finally, reading off all coefficients of $M(d)$ with $c = 0$
yields the complete list of ``unmixed'' two-point invariants,
which in turn yield the structure constants of thesmall quantum
product (via the associativity). 

\medskip

This algorithm is very easy to implement. For example, when $X$
is a quintic hypersurface in ${\bf P}^6$, it gives the following
products:

\begin{equation*}
  \begin{matrix}
  H*H   & = & H^2 & + & 120q     &   &                 &   &\\
  H*H^2 & = & H^3 & + & 770qH    &   &                 &   &\\
  H*H^3 & = & H^4 & + & 1345qH^2 & + & 211,200q^2      &   &\\
  H*H^4 & = & H^5 & + & 770qH^3  & + & 692,\!500q^2H   &   &\\
  H*H^5 & = &     &   & 120qH^4  & + & 211,\!200q^2H^2 & + & 31,\!320,\!000q^3
  \end{matrix}
\end{equation*}

As a typical application, note that the last number implies the
following interesting bit of enumerative geometry:

\begin{cor} The expected number of twisted cubics through two general
points of of a quintic five-fold $X\subset\P[6]$ is:
\[
  2,\!088,\!000.
\]
\end{cor}

One similarly may produce the expected numbers of 
rational normal curves of degree $d$ passing through $2$ general 
points of a hypersurface of degree $2d-1$ in $\P[2d]$ for 
any~$d$.

\providecommand{\bysame}{\leavevmode\hbox to3em{\hrulefill}\thinspace}

\end{document}